\documentclass[12pt]{amsart}
\usepackage{amsmath,amssymb,amscd}
\usepackage[colorlinks=true, pdfstartview=FitV, linkcolor=blue, citecolor=blue, urlcolor=blue]{hyperref}
\usepackage{mathpazo}

\newcommand{\treg}{\text{reg}}

\theoremstyle{plain}

\newtheorem*{nonumtheorem}{Theorem}

\newtheorem*{nonumcorollary}{Corollary}

\theoremstyle{definition}

\begin{document}
\title[A number field example]{On Zagier's conjecture for $L(E,2)$:\\a number field example}
\author{Jeffrey Stopple}
\begin{abstract}
We work out an example, for a CM elliptic curve $E$ defined over a real quadratic field $F$, of Zagier's conjecture.  This relates $L(E,2)$ to values of the elliptic dilogarithm function at a divisor in the Jacobian of $E$ which arises from $K$-theory.
\end{abstract}
\address{Mathematics Department, UC Santa Barbara, Santa Barbara CA 93106}
\email{stopple@math.ucsb.edu}
\maketitle

\subsection*{Introduction}
Recall
that the classical Euler dilogarithm is defined by
\begin{align*}
Li_2(z)=&\sum_{n=1}^\infty \frac{z^n}{n^2}\qquad |z|<1\\
=&\int_0^z-\frac{\log(1-t)}{t} dt\qquad z\in\mathbb C\backslash[1,\infty)
\end{align*}
after analytic continuation.  The Bloch-Wigner dilogarithm
\[
D(z)=\text{Im}(Li_2(z))+\log|z|\arg(1-z)
\]
is well defined independent of path used to continue $Li_2$ and $\arg$.  For a torus
$\mathbb C/\Lambda$, $\Lambda=[\omega_1,\omega_2]$ corresponding to a point $\tau$ in $\mathcal H$,
we have the
$q$-symmetrized, or elliptic, dilogarithm
\[
D_q(z)=\sum_{n\in\mathbb Z} D(zq^n)\quad q=\exp(2\pi i\tau)\quad z\in \mathbb C^\times/q^{\mathbb
Z}\cong
\mathbb C/\Lambda.
\]

In their paper on Zagier's Conjecture, Goncharov and Levin prove the following theorem  \cite[Theorem 1.1]{GL} about the value at $s=2$ of the $L$-function of an elliptic curve:
\begin{nonumtheorem}  Let $E$ be a (modular) elliptic curve over $\mathbb Q$.  Then there exists a $\mathbb Q$-rational divisor $\mathcal P=\sum a_i(P_i)$ satisfying
\begin{equation}
\sum a_i P_i\otimes P_i\otimes P_i=0\qquad \text{ in }S^3 J(E),\tag{a}\label{Taga}
\end{equation}
For any valuation $v$ of the field $\mathbb Q(\mathcal P)$, and $h_v$ the corresponding canonical height, 
\begin{equation}
\sum a_i h_v(P_i)\cdot P_i=0\qquad\text{ in }J(E)\otimes \mathbb R,\tag{b}\label{Tagb}
\end{equation}
as well as a certain third condition (c) at primes where $E$ has split multiplicative reduction.  For such a divisor
\[
L(E,2)\sim_{\mathbb Q^\times}\pi \mathcal D_q(\mathcal  P).
\]
\end{nonumtheorem}
\noindent
The authors comment 
\begin{quote}\emph{
\lq\lq The conditions (a) and (b) were guessed by Zagier several years ago after
studying the results of the computer experiments with $\mathbb Q$-rational points on
some elliptic curves, which he did with H. Cohen.\rq\rq}
\end{quote}
As a consequence, they deduce \cite[Corollary 1.3]{GL}
\begin{nonumcorollary}  Let $E$ be an elliptic curve over $\mathbb Q$.  Let us assume that the image of $K_2(E)_{\mathbb Z}\otimes\mathbb Q$ under the regulator map  is $L(E,2)\cdot \mathbb Q$.  (This is part of the Bloch-Beilinson conjecture.)  Then for any $\mathbb Q$-rational divisor $\mathcal P$ on $E(\bar{\mathbb Q})$ satisfying the conditions (\ref{Taga}), (\ref{Tagb}), and (c) above, one has 
\[
r\cdot  L(E,2)=\pi D_q(\mathcal P)
\]
where  $r$ is a rational number, perhaps equal to $0$.
\end{nonumcorollary}
\noindent 
They remark 
\begin{quote}\emph{\lq\lq Theorem 1.1 and Corollary 1.3 have analogs for an elliptic curve
over any number field. Its formulation is an easy exercise to the reader.\rq\rq}\end{quote}

In \cite[\S 1.3]{GL} they work out an example for the elliptic curve given by $y^2-y=x^3-x$.  Nonetheless, examples in this subject are scarce and the theory is more than a little intimidating.  The purpose of this note is to work out an example for a curve over a number field, following the philosophy of \cite{Dinakar}, 
\begin{quote}\emph{\lq\lq In
general, the more concrete one is able to make the [Borel] regulator map, the more explicit the information one is able to
extract from it.\rq\rq}
\end{quote}

All of the calculations were done with PARI.

\subsection*{Notation}
Let $F$ be the field $\mathbb Q(\sqrt{5})$, and
$w=\frac{1+\sqrt{5}}{2}.$  Consider the elliptic curve $E$ defined over $F$:
\[E:\,\,y^2+y=x^3+w\, x^2-(93+163w)\,x+(669+1076w).
\] The discriminant is
$-5^37^3=-42875$, and the $j$-invariant is 
\[
j(E)=-32604160-52756480w=(-128-224w)^3.
\]
This is equal to
$j(\mathcal O_K)$, where $K=\mathbb Q(\sqrt{-35}).$  Thus  the curve has complex multiplication by
the ring of integers
$\mathcal O_K$.  The two embeddings of $F$ into $\mathbb R$ give two lattices 
\[
\Lambda=[1,\tau],\quad 35\tau^2+35\tau+9=0,\qquad
\Lambda^\prime=[1,\tau^\prime],\quad 7{\tau^\prime}^2+7\tau^\prime+3=0.
\] The fact that
$\tau^\prime=5\tau+2$ show that $E$ is isogenous to its Galois conjugate, so it is a $\mathbb Q$-curve in the sense of \cite{Gr}.  $K$ has class number 2 since the $j$-invariant is quadratic,
while $F$ has class number 1.  Note that since $E$ has complex multiplication, it has only additive bad
reduction and we can ignore condition (c) in the theorem.
\subsection*{$L$-function computation}  The curve $E$ is in fact the canonical $\mathbb Q$-curve
(Theorem.11.2.4 of \cite{Gr}) for this discriminant, which is convenient for calculating values of the
$L$-function.   The Hecke character $\psi$ on the Hilbert class field
$H$ factors through norms from $H$ to $K$.  The Euler product at $s=2$ converges too slowly to
be of use.  So we use the functional equation to convert the value at $s=2$ to the leading Taylor coefficient at $s=0$.  Since the field is quadratic there is a second order zero at $0$.  Thus we are computing the value $L(E,0)^{(2)}$, or up to appropriate powers of $\pi$ and rational multiples, the value of the \lq completed\rq\ $\Lambda(E,s)$ at $s=0$.   

Following ideas of Cremona \cite{CW} we write the
$L$ function as the Mellin transform of a Maass form on $\mathcal H^3$, with a Fourier series
involving $K$-Bessel functions.  
Although $K$ has class number 2, the Maass form is a \lq CM\rq
form, so its Fourier coefficients are supported on the principal ideal class.  We split the
integral at the symmetry point, use the functional equation, and integrate by parts.  To get 28 digits of accuracy we computed the Dirichlet
series coefficients for primes less than 30,000.  The values of
$\Lambda(E,s)$ at $s=0$ require evaluating, for 30,000 different $x$ values,
\[
\int_x^\infty K_0(t)/t \,dt.
\] For $x\leq 3$ or $15\leq x$, we can take an asymptotic expansion for $K_0(t)$ and integrate term
by term to get an asymptotic expansion for the function.  For $3<x<15,$ we need to numerically
integrate from
$x$ to the next integer ceil$(x)$, and use a table lookup for
$\int_{\text{ceil}(x)}^\infty K_0(t)/t \,dt.$ Eventually we find
\[ L(E,0)^{(2)}=691.9884130215329129184499757.
\]

\subsection*{Regulator computation}
Let $P=[7+9w,17+35w]$ and let
$Q=[12-w,32-20w]$.  These points seem to generate the free part of the group $E(F).$  The curve 
has a large number of integral points, (1)-(14) in Table \ref{Ta:Points}.  In order to find solutions $a_i$ to the equations (\ref{Taga}), (\ref{Tagb}) in the construction
of Goncharov and Levin,  one needs a relatively large number of points whose local heights are supported on a
relatively small number of primes.  We consider also the points (15)-(22) in Table
\ref{Ta:Points}.  The local nonarchimedean height functions are supported on the primes
$2,\sqrt{5},7,\pi_{11}$ and $\pi_{59}$, where $\pi_{11}$ and $\pi_{59}$ are primes above 11 and 59
in $\mathcal O_F.$

\begin{table}
\begin{flushleft}
\begin{tabular}{l r l} (1)&$[0]P+[1]Q=$&$[12-w,32-20w]$\\ (2)&$[1]P+[0]Q=$&$[7+9w,17+35w]$\\
(3)&$[0]P+[2]Q = $&$[-4 - 11w, 11 + 8w]$\\ (4)&$[1]P + [1]Q = $&$[7 + 2w, -11 + 7w]$\\ (5)&$[2]P
+[0]Q = $&$[3 + 5w, 2 + w]$\\ (6)&$[1]P + [2]Q = $&$[42 - 26w, -333 + 175w]$\\ (7)&$[2]P + [1]Q =
$&$[2 + 4w, 2 + 5w]$\\ (8)&$[2]P + [2]Q = $&$[3 + 4w, -4 - w]$\\ (9)&$[3]P - [1]Q = $&$[1624 -
957w, -75625 + 46340w]$\\ (10)&$[3]P + [1]Q = $&$[-5w, 24 + 28w]$\\ (11)&$[4]P + [1]Q = $&$[27 -
26w, -223 + 95w]$\\ (12)&$[4]P + [2]Q = $&$[46 - 22w, 331 - 205w]$\\ (13)&$[4]P + [3]Q = $&$[67 +
99w, 957 + 1525w]$\\ (14)&$[5]P + [4]Q = $&$[250362 - 154726w, -147263008 + 91013545w]$\\ 
(15)&$[1]P - [1]Q =$&$ [(14 + 24w)/5, ...]$\\ (16)&$[1]P - [2]Q =$&$ [(2527 + 6584w)/3481, ...]$\\
(17)&$[3]P +[0]Q=$&$ [(217 - 31w)/16, ...]$\\ (18)&$[3]P + [2]Q =$&$ [(392 + 529w)/121, ...]$\\
(19)&$[4]P + [4]Q =$&$ [(13627 + 13872w)/3481, ...]$\\ (20)&$[5]P + [2]Q =$&$ [(17367 +
12464w)/3481, ...]$\\ (21)&$[6]P +[0]Q= $&$[(792753 + 52969w)/222784,...]$\\ (22)&$[6]P + [4]Q =$&$
[(1700 + 1357w)/605, ...]$\\ 
\end{tabular}
\caption{}\label{Ta:Points}
\end{flushleft}
\end{table}
Since $E$ has rank (at least) 2, it will be convenient to revise our notation for a divisor
\[
\mathcal P=\sum_{k,l} a_{k,l}\left([k]\cdot P+[l]\cdot Q\right),
\]
where $k$ and $l$ are restricted to the values in Table \ref{Ta:Points}.  The condition (\ref{Taga}) becomes the four equations
\[
\sum a_{k,l}\cdot k^3=0,\quad
\sum a_{k,l}\cdot k^2l=0,\quad
\sum a_{k,l}\cdot kl^2=0,\quad
\sum a_{k,l}\cdot l^3=0.
\]
Meanwhile condition (\ref{Tagb}) becomes, for the nonarchimedean heights,
\begin{gather*}
\sum a_{k,l} h_v([k]P+[l]Q)\cdot k=0\\
\sum a_{k,l} h_v([k]P+[l]Q) \cdot l=0,
\end{gather*}
ten more equations as $v$ ranges over the five primes $2,\sqrt{5},7,\pi_{11}$ and $\pi_{59}$.
To compute the height functions, we used \cite{AECII} and the reference therein, particularly \cite{Silverman}.  These equations are defined over $\mathbb Z$, so we get integral solutions.  Surprisingly,
the solution space is 10, not 8 dimensional.  

\begin{table}
\begin{center}
\begin{tabular}{ c c c c c c c c } 3& 7& -8& -5& -10& -7& -8& 7\\ 2& -23& -11& -45& -48& 18& -181&
-33\\ 1& -1& -9& -4& -8& 1& -33& 3\\ 6& -1& -11& 15& -30& 1& -10& 13\\ 3& 5& -12& -17& -13& -1&
-45& -7\\ -2& -1& 16& 3& -4& -3& 37& 0\\ -2& -4& 0& 14& -14& 4& -14& 26\\ -3& -1& 18& 1& 26& 7& 52&
-22\\ 0& 0& 1& 0& -2& 0& -1& 1\\ -2& 1& 0& 1& 0& -1& 1& 0\\ 0& 0& -1& -2& -1& 0& 5& -2\\ 0& -1& 1&
0& 0& -3& -2& -1\\ 1& 1& -1& 1& -1& -1& 1& -1\\ 0& 0& -1& 0& 1& 1& -3& 0\\ 0& -2& 0& 2& 0& 2& 2& 0\\
0& 0& -2& -3& -5& -1& -8& 0\\ 0& 0& 0& -6& -9& 0& -9& -3\\ 0& 0& 0& 0& 0& 0& 0& -4\\ 0& 0& -2& -1&
-2& -1& -5& 1\\ 0& 0& 2& -1& -1&1& 2& -2\\ 0& 0& 0& 2& 3& 0& 3& 1\\ 0& 0& 0& 0& 0& 0& 0& 2\\ 
\end{tabular}
\caption{8 divisors supported on the 22 points}\label{Ta:Symbol}
\end{center}
\end{table}

\begin{table}
\begin{center}
\begin{tabular}{ r r }
$1.7 \times 10^{-100} $&
$-2.9 \times 10^{-105}$\\
$3.657296793764310936796018961 $&
$-4.861051673717091496858129462$\\
$-3.657296793764310936796018961$ &
$4.861051673717091496858129462$\\
$25.64710971025614581418182019$ &
$1.387883495576657586500860340$\\
$3.657296793764310936796018961$ &
$-4.861051673717091496858129462$\\
$-3.657296793764310936796018961$ &
$4.861051673717091496858129462$\\
$35.41524521159629806450776657$ &
$0.2301607695298462830484372999$\\
$29.30440650402045675097783915$ &
$-3.473168178140433910357269121$\\ 
\end{tabular}
\caption{Regulators in $\mathbb R^2$}\label{Ta:Reg}
\end{center}
\end{table}

In this solution space we next seek integral solutions to the equations over $\mathbb R$
\begin{gather*}
\sum a_{k,l} h_\infty([k]P+[l]Q)\cdot k=0\\
\sum a_{k,l} h_\infty([k]P+[l]Q) \cdot l=0,
\end{gather*}
 for each of the two archimedean valuations.  Since the condition (\ref{Taga}) kills the global
height, it suffices to find solution for just one infinite prime, and we can use the other one as a
check that our calculation is correct.  Finding integral solutions to equations given by real
(floating point) numbers is tricky.  The easiest way seems to be (following Zagier \cite{Z}) to use
the LLL algorithm.  We get an 8 dimensional solution space of integral vectors, see Table
\ref{Ta:Symbol}.
Corresponding to column
$j$ in Table \ref{Ta:Symbol} above is the divisor 
\[
\mathcal P(j)=\sum_{k,l} a_{k,l,j}([k]P+[l]Q).
\]

For each point $[k]P+[l]Q$ we compute $z_{k,l}$ modulo $\Lambda$ and $z_{k,l}^\prime$ modulo $\Lambda^\prime$. 
Let $q=\exp(2\pi i\tau)$,
$q^\prime=\exp(2\pi i\tau^\prime)$.  
Then corresponding to column
$j$ in Table \ref{Ta:Symbol} we compute the vector in $\mathbb R^2$ given by:
\begin{multline*}
\{\treg(\mathcal P(j)),\treg(\mathcal P(j))^\prime\}=\\
\sum_{k,l} a_{k,l,j} \{D_q(\exp(2\pi i z_{k,l})/\pi,D_{q^\prime}(\exp(2\pi i
z^\prime_{k,l})/\pi\}.
\end{multline*}
Working with 100 digits, (displaying 28), we get the row vectors in Table \ref{Ta:Reg}.

\subsection*{Comparison}
With 8 regulator vectors in $\mathbb R^2$, there are, up to sign, 28 choices for a $2\times 2$
determinant $R_{m,n}$ of the rows $m$ and $n$, of which 13 visibly have determinant
equal 0.  For the remaining 15 pairs we get that $R_{m,n}/L(E,0)^{(2)}$ appears to be rational. 
The pair
$(4,7)$ gives
\[ -0.06250000000000000000000000268\approx-\frac{1}{16}
\] while $(2,7)$, $(3,7)$, $(5,7)$, $(6,7)$ all give plus or minus
\[ 0.2500000000000000000000000107 \approx\frac{1}{4}
\] and $(2,4)$, $(2,8)$, $(3,4)$, $(3,8)$, $(4,5)$, $(4,6)$, $(4,8)$, $(5,8)$, $(6,8)$, and $(7,8)$
all give plus or minus
\[ 0.1875000000000000000000000080\approx\frac{3}{16}.
\]

The close agreement with a rational number of small denominator serves as confirmation the calculations are correct.

\end{document}